%
\input amstex.tex
\loadmsam
\loadmsbm
\loadbold
\input amssym.tex
\baselineskip=13pt plus 2pt
\documentstyle{amsppt}
\pageheight{45pc}
\pagewidth{33pc}

\magnification=1200
\overfullrule=0pt

\def\enddemos{\hfill$\square$\enddemo}
\def\hb{\hfil\break}
\def\n{\noindent}
\def\ov{\overline}
\def\smatrix{\smallmatrix}

\def\pmatrix{\left(\smatrix}
\def\endpmatrix{\endsmallmatrix\right)}
\def\liminj{\lim\limits_{\longrightarrow}}
\def\limproj{\lim\limits_{\longleftarrow}}

\def\upi{\pmb{\pi}}

\def\upsi{{\pmb{\psi}}}
\def\bs{\backslash}

\def\f{\varphi}

\def\wh{\widehat}
\def\A{\Bbb A}
\def\C{\Bbb C}

\def\1{\bold 1}

\def\ch{\operatorname{ch}}
\def\Gal{\operatorname{Gal}}

\def\GL{\operatorname{GL}}
\def\U{\operatorname{U}}

\def\Hom{\operatorname{Hom}}

\def\PGL{\operatorname{PGL}}

\def\mod{\operatorname{mod}}
\def\ind{\operatorname{ind}}
\def\Ind{\operatorname{Ind}}
\def\tr{\operatorname{tr}}

\def\LieG{\Cal G}
\def\LieH{\Cal H}

\def\LieM{\Cal M}

\def\diag{\operatorname{diag}}

\leftheadtext{Yuval Z. Flicker}
\rightheadtext{Characters, genericity, and multiplicity one for U(3)}
\topmatter
\title\nofrills{CHARACTERS, GENERICITY, AND MULTIPLICITY ONE FOR U(3)}
\endtitle
\author Yuval Z. Flicker\endauthor
\footnote"~"{\n Department of Mathematics, The Ohio State University,
231 W. 18th Ave., Columbus, OH 43210-1174.\hb
Email: flicker\@math.ohio-state.edu\hb
\indent Partially supported by a Lady Davis Visiting Professorship at 
the Hebrew University and Max-Planck scholarship at MPI, Bonn.$\qquad$ 
I am grateful to J. Bernstein and G. Savin for constructive criticism.\hb 
}

\abstract Let $\upsi:U\to\C^\times$ 
be a generic character of the unipotent 
radical $U$ of a Borel subgroup of a quasisplit $p$-adic group $G$.
The number (0 or 1) of $\upsi$-Whittaker models on an admissible 
irreducible representation $\pi$ of $G$ was expressed by Rodier in terms 
of the limit of values of the trace of $\pi$ at certain measures
concentrated near the origin. An analogous statement holds in the twisted 
case. This twisted analogue is used in [F], p. 47, to provide a local proof 
of the multiplicity one theorem for U(3). This asserts that each discrete 
spectrum automorphic representation of the quasisplit unitary group U(3) 
associated with a quadratic extension $E/F$ of number fields occurs in the 
discrete spectrum with multiplicity one. It is pointed out in [F], p. 47, 
that a proof of the twisted analogue of Rodier's theorem does not appear 
in print. It is then given below. Detailing this proof is necessitated in
particular by the fact that the attempt in [F], p. 48, at a global proof 
of the multiplicity one theorem for U(3), although widely quoted, is 
incomplete, as we point out here.
\endabstract
\endtopmatter

\document

Let $E/F$ be a quadratic extension of $p$-adic fields, $p\not=2$,
and consider the basechange lifting from the quasisplit unitary
group $H=\U(3,E/F)$ in 3 variables to $G=\GL(3,E)$.
Our main result is that there exists a family of suitably related
functions $\psi'_{H,n}$ on $H=\U(3,E/F)$ and $\psi'_n$ on
$G=\GL(3,E)$ which are supported near the origin, such that the traces
$\tr\pi_H(\psi'_{H,n}dh)$ and twisted traces
$\tr\pi(\psi'_ndg\times\sigma)$ stabilize
for sufficiently large $n$ and become equal to the Whittaker multiplicity
(multiplicity in the space of Whittaker vectors) of the irreducible
admissible representations $\pi_H$ of $H$ and $\pi$ of $G$. This is a 
twisted analogue of a theorem of Rodier for the involution that defines 
$H$ in $G$.

Our motivation for considering such a twisted analogue of Rodier's theorem
is that we found a gap in an attempted global proof (of [F], p. 48) of 
multiplicity one theorem for automorphic representations of $\U(3,E/F)(\A_F)$.
We use our result to provide a local proof, as suggested in [F], p. 47.

This gap is different than the observation of Harder ([H], p. 173).
The latter deals with the question of which representation in a packet
is automorphic, assuming multiplicity one holds.

In Section 1 we state a theorem of Rodier and explain what are Whittaker 
models, characters, and the measures concentrated near the origin which 
enter the statement. We then state its twisted analogue of interest to us.
In Section 2 we use these theorems to complete a key step in the proof of
[F], p. 47, of multiplicity one theorem for the global quasisplit unitary 
group U(3). In Section 3 we explain why the global proof of [F], p. 48, of 
this multiplicity one theorem is not complete. In Section 4 we recall the 
main lines of Rodier's proof. In Section 5 the twisted analogue is reduced 
to Rodier's theorem. This completes the proof of the theorems of Sections 1 
and 2. In the Appendix, Section 6, we derive another description of the
Whittaker multiplicity, in terms of the coefficient of the
regular orbit in the germ expansion of the character.

\heading{1. Whittaker models and characters}\endheading

Rodier's theorem [R], p. 161, (for a split group $H$) computes the number 
of $\upsi_H$-Whittaker models of the admissible irreducible representation 
$\pi_H$ of $H$ in terms of values of the character $\tr\pi_H$ or 
$\chi_{\pi_H}$ of $\pi_H$ at the measures $\psi_{H,n}dh$ which are 
supported near the origin.

We proceed to explain the notations to be used in Rodier's theorem. For 
simplicity and clarity, instead of working with a general connected reductive 
(quasi) split $p$-adic group $H$ as in [R], we let $H$ be a specific unitary
group. To define it, we take $G=\GL(r,E)$, where $E/F$ is a quadratic extension
of $p$-adic fields of characteristic zero, $p\not=2$. Let $x\mapsto\ov{x}$ 
denote the generator of $\Gal(E/F)$. For $g=(g_{ij})$ in $G$ we put 
$\ov{g}=(\ov{g}_{ij})$ and ${}^tg=(g_{ji})$. Then $\sigma(g)
=J^{-1}{}^t\ov{g}^{-1}{}J$, $J=((-1)^{i-1}\delta_{i,r+1-j})$, defines an
involution $\sigma$ on $G$. The group $H=G^\sigma$ of $g\in G$ fixed by 
$\sigma$ is a quasisplit unitary group $\U(r,E/F)$. 

Denote by $\upsi_H:U_H\to\C^\times$ a character on the unipotent upper 
triangular subgroup $U_H$ of $H$. It is necessarily unitary, namely its
values lie in the unit circle $\C^1=\{z\in\C;|z|=1\}$.

We assume that $\upsi_H$ is generic, namely nontrivial on each simple root 
subgroup. There is only one orbit of generic $\upsi_H$ under the action 
of the diagonal subgroup of $H$ on $U_H$ by conjugation. Hence we can
and will work with the specific character $\upsi_H:U_H\to\C^1$, defined 
by $\upsi_H((u_{ij}))=\psi(\sum_{1\le j<r}u_{j,j+1})$. Here
$\psi:F\to\C^1$ is an additive character which is 1 on $R$ and
not identically 1 on $\upi^{-1}R$. Further, $R$ is the ring of integers 
of $F$, and $\upi$ is a generator of the maximal ideal of $R$.
Note that $u_{r-j,r-j+1}=\ov u_{j,j+1}$.

By $\upsi_H$-{\it Whittaker vectors} we mean vectors in the space of 
the induced representation $\ind_{U_H}^H(\upsi_H)$. They are the functions 
$\f:H\to\C$ with $\f(uhk)=\upsi_H(u)\f(h)$, $u\in U_H$, $h\in H$, 
$k\in K_{\f}$, where $K_{\f}$ is a compact open subgroup of $H$ 
depending on $\f$, which are compactly supported on $U_H\bs H$. 
The group $H$ acts by right translation. 

The multiplicity $\dim_\C\Hom_H(\ind_{U_H}^H\upsi_H,\pi_H)$ 
of any irreducible admissible 
representation $\pi_H$ of $H$ in the space of $\upsi_H$-Whittaker vectors 
is known to be 0 or 1. In the latter case we say that $\pi_H$ has a 
$\upsi_H$-Whittaker model or that it is $\upsi_H$-{\it generic}.

Let $\LieG_0$ be the ring of $r\times r$ matrices with entries in the ring 
of integers $R_E$ of $E$. It is a subring of the ring $\LieG$ of $r\times r$
matrices with entries in $E$. Let $d\sigma$ be the involution 
$d\sigma(X)=-J^{-1}{}^t\ov XJ$. Its set of fixed points in $\LieG$ is denoted 
by $\LieH$. Let $\LieH_0=\LieH\cap\LieG_0$ be the set of $X$ in $\LieG_0$ 
fixed by the involution $d\sigma$.
Write $H_n=\exp(\LieH_n)$, $\LieH_n=\upi^n\LieH_0$. For $n\ge 1$ we have 
$H_n={}^tU_{H,n}A_{H,n}U_{H,n}$, where $U_{H,n}=U_H\cap H_n$, and $A_{H,n}$ 
is the group of diagonal matrices in $H_n$. Define a character 
$\psi_{H,n}:H\to\C^1$, supported on $H_n$, by
$\psi_{H,n}({}^tbu)=\psi(\sum_{1\le j< r}u_{j,j+1}\upi^{-2n})$,
at ${}^tb\in{}^tU_{H,n}A_{H,n}$, $u=(u_{ij})\in U_{H,n}$.
Alternatively, by
$$\psi_{H,n}(\exp X)=\ch_{\LieH_n}(X)\psi(\tr[X\upi^{-2n}\beta_H]),$$
where $\ch_{\LieH_n}$ indicates the characteristic function of $\LieH_n
=\upi^n\LieH_0$ in $\LieH$, and $\beta_H$ is the $r\times r$ matrix whose 
nonzero entries are 1 at the places $(j,j-1)$, $1<j\le r$. Denote by
$e_{H_n}$ the constant measure of volume one supported on the compact
subgroup $H_n$ in the Hecke algebra of $H$, i.e., 
$e_{H_n}=|H_n|^{-1}\ch_{H_n}dh$, where $|H_n|$ denotes the volume of
$H_n$ in $dh$. 

Since $\pi_H$ is admissible, for each test measure $f_H\,dh$ 
($dh$ is a Haar measure on $H$ and $f_H$ is a locally constant compactly 
supported complex valued function on $H$), the image of $\pi_H(f_H\,dh)$ is 
finite dimensional and its trace $\tr\pi_H(f_H\,dh)$ is finite. 

Rodier's theorem asserts

\proclaim{Theorem 1} The multiplicity 
$\dim_{\C}\Hom_H(\ind^H_{U_H}\upsi_H,\pi_H)$
is equal to
$$\lim_n\tr\pi_H(\psi_{H,n}e_{H_n}).$$
\endproclaim

In fact the limit stabilizes for sufficiently large $n$. Throughout this paper
``= $\lim_na_n$'' means ``equals $a_n$ for all sufficiently large $n$''.
\medskip

We need a twisted analogue of Rodier's theorem. It can be described as follows.

Let $\pi$ be an admissible irreducible representation of $G$ which is 
$\sigma$-invariant: $\pi\simeq{}^\sigma\pi$, where 
${}^\sigma\pi(\sigma(g))=\pi(g)$.
Then there exists an intertwining operator $A:\pi\to {}^\sigma\pi$, with 
$A\pi(g)=\pi(\sigma(g))A$ for all $g\in G$. Since $\pi$ is irreducible, by 
Schur's lemma $A^2$ is a scalar which we may normalize by $A^2=1$. Thus $A$ 
is unique up to a sign. Denote by $G'$ the semidirect product 
$G\rtimes\langle\sigma\rangle$. Then $\pi$ extends to $G'$ by $\pi(\sigma)=A$.

Define $\psi_E:E\to\C^1$ by $\psi_E(x)=\psi(x+\ov x)$. Define a character 
$\upsi:U\to \C^1$ on the unipotent upper triangular subgroup $U$ of $G$ by 
$\upsi((u_{ij}))=\psi_E(\sum_{1\le j<r}u_{j,j+1}).$
This one dimensional representation has the property that  
$\upsi(\sigma(u))=\upsi(u)$ for all $u$ in $U$. Note that 
$\upsi(u)=\upsi_H(u^2)$ at $u\in U_H=U\cap H$. There is only one orbit of 
generic $\sigma$-invariant characters on $U$ under the adjoint action of the 
group of $\sigma$-invariant diagonal elements in $G$.

Suppose that $\pi$ is $\upsi$-generic, namely 
$\Hom_{\C}(\ind_U^G\upsi,\pi)\not=\{0\}$.
Here $\ind_U^G\upsi$ consists of the functions $\f:G\to\C$ with 
$\f(ug)=\upsi(u)\f(g)$, $u\in U$, $g\in G$, which are compactly supported
on $U\bs G$. Then we normalize $A$ by $A\f'={}^\sigma \f'$, where 
${}^\sigma \f(g)=\f(\sigma(g))$, on the image $\f'$ in $\pi$ of the $\f$.

Write $G_n=\exp(\LieG_n)$, where $\LieG_n=\upi^n\LieG_0$. For $n\ge 1$ we 
have $G_n={}^tU_nA_nU_n$, where $U_n=U\cap G_n$, and $A_n$ is the group of 
diagonal matrices in $G_n$. Define a character $\psi_n:G\to\C^1$
supported on $G_n$ by 
$\psi_n({}^tbu)=\psi_E(\sum_{1\le j<r}u_{j,j+1}\upi^{-2n})$
where ${}^tb\in{}^tU_nA_n$, $u=(u_{ij})\in U_n$.
Alternatively, $\psi_n:G\to\C^1$ is defined by
$$\psi_n(\exp X)=\ch_{\LieG_n}(X)\psi_E(\tr[X\upi^{-2n}\beta])$$
where $\beta$ is the $r\times r$ matrix with entries 1 at the
places $(j,j-1)$, $1<j\le r$, and 0 elsewhere.

A first version of a $\sigma$-twisted analogue of Rodier's theorem
asserts that the Whittaker multiplicity is equal to the twisted trace
$\tr\pi(\psi_ne_{G_n}\times\sigma)$ for all sufficiently large $n$,
where $e_{G_n}$ is defined analogously to $e_{H_n}$. 

A more useful version for us is stated in terms of the twisted character 
$\chi_\pi^\sigma$. Let us first restate Rodier's theorem this way.

Denote by $\chi_{\pi_H}$ the {\it character} of $\pi_H$. It is a complex 
valued conjugacy invariant function on $H$ which is locally constant on 
the regular set and locally integrable on $H$ (Harish-Chandra [HC], Theorem 
1) defined by $\tr\pi_H(f_H\,dh)=\int_H\chi_{\pi_H}(h)f_H(h)dh$ for all 
$f_H\,dh$.

Then Rodier's theorem can be stated as asserting that the Whittaker 
multiplicity of Theorem 1 is equal to
$$\lim_n\int_H\chi_{\pi_H}(h)\psi_{H,n}(h)e_{H_n}(h).$$

Analogously, the twisted character $\chi_\pi^\sigma$ of $\pi$ is a complex 
valued $\sigma$-conjugacy invariant function on $G$ (that is, its value 
on $\{hg\sigma(h)^{-1}\}$ is independent of $h\in G$) which is locally 
constant on the $\sigma$-regular set ($g$ with regular $g\sigma(g)$), 
locally integrable (Clozel [C], Thm 1, p. 153) and defined by 
$\tr\pi(f\,dg)A=\int_G\chi_\pi^\sigma(g)f(g)dg$ for all test 
measures $f\,dg$.

The $\sigma$-twisted analogue of Rodier's theorem of interest to us is 
as follows. Let $e_{G_n^\sigma}$ denote the constant measure of volume 1 
supported on the compact subgroup $G_n^\sigma=\{g=\sigma g;g\in G_n\}$ of 
$G$. As $G_n^\sigma$ is $H_n$, $e_{G_n^\sigma}$ lies in the Hecke algebra
 of $H$.

\proclaim{Theorem 2} The multiplicity
$\dim_\C\Hom_{G'}(\ind_U^G\upsi,\pi)=\dim_\C\Hom_{G}(\ind_U^G\upsi,\pi)$
is equal to
$\int_{G_n^\sigma}\chi_\pi^\sigma(g)\psi_n(g)e_{G^\sigma_n}(g)$
for all sufficiently large $n$.
\endproclaim

\remark{Remark} Recall that 
$\Hom_G(\pi_1,\pi_2^\vee)=\Hom_G(\pi_2,\pi_1^\vee)$ 
as both spaces can be identified with the space of $(\pi_1,\pi_2)$-invariant
bilinear forms. The contragredient $(\ind_H^G\rho)^\vee$ is 
$\Ind_H^G({{\Delta_G}\over{\Delta_H}}\rho^\vee)$ ([BZ1], 2.25(c); Ind 
indicates non compact induction, $H$ is a closed subgroup of $G$, $\rho$
is a representation of $H$, $\pi$ of $G$). Frobenius reciprocity asserts
$$\Hom_G(\ind_H^G\rho,\pi)=\Hom_H({{\Delta_H}\over{\Delta_G}}\rho,\pi|H)$$
and equivalently $\Hom_G(\pi,\Ind_H^G\rho)=\Hom_H(\pi|H,\rho)$ ([BZ1], (2.28)
and (2.29)).\endremark

\heading{2. Multiplicity one for U(3)}\endheading

Let us recall how the Theorems are used in the proof of Proposition 3.5 
of [F], p. 47. Thus in this Section we work only with $H=\U(r,E/F)$ and 
$G=\GL(r,E)$, $r=2$ or $3$, of Section 1. We are given a square integrable 
irreducible admissible representation $\rho$ of the quasisplit group 
$\U(2,E/F)$. Its stable basechange to $\GL(2,E)$ is denoted by $\tau$.
The unstable basechange is $\tau\otimes\kappa$. Let $\pi=I(\tau\otimes\kappa)$
be the normalizedly induced representation of $\GL(3,E)$.
This $\pi$ is invariant under the involution $\sigma$ (ours and of [F]). 
It is generic. For all matching measures $fdg$ and $f_Hdh$ on $G=\GL(3,E)$ 
and $H=\U(3,E/F)$, using an identity of trace formulae and orthogonality 
relations for characters [F] obtains an identity 
$$\tr\pi(fdg\times\sigma)=(2m+1)\sum_{\pi_H}\tr\pi_H(f_Hdh).$$ 
The sum ranges over finitely many (in fact, two times the cardinality of 
the packet of $\rho$) inequivalent square integrable irreducible admissible 
representations $\pi_H$ of $\U(3,E/F)$. The number $m$ is a nonnegative 
integer, independent of $\pi_H$. 

\proclaim{Proposition 3.5 of [F]} The nonnegative integer $m$ is zero, and 
there is a unique generic $\pi_H$ in the sum. The other $2[\{\rho\}]-1$ 
representations $\pi_H$ are not generic. 
\endproclaim

Note that our $\pi$, $\pi_H$, $fdg$, $f_Hdh$, $G$, $H$, are denoted in [F] by 
$\Pi$, $\pi$, $\phi dg'$, $fdg$, $G'$, $G$.

\demo{Proof} The identity for all matching test measures $fdg$ and $f_Hdh$ 
implies an identity of characters: 
$$\chi_\pi^\sigma(\delta)=(2m+1)\sum_{\pi_H}\chi_{\pi_H}(\gamma)$$
for all $\delta\in G=\GL(3,E)$ with regular norm $\gamma\in H=\U(3,E/F)$.
Note that $\delta\mapsto\chi_\pi^\sigma(\delta)$ is a stable 
$\sigma$-conjugacy class function on $G$, while 
$\gamma\mapsto\sum_{\pi_H}\chi_{\pi_H}(\gamma)$ is a stable conjugacy class 
function on $H$.
We use Theorem 2 with $G=\GL(3,E)$ and $H=G^\sigma$. Then $G_n^\sigma=H_n$.
On $\delta\in G_n^\sigma$, the norm $N\delta$ of the stable $\sigma$-conjugacy 
class $\delta$ is just the stable conjugacy class of $\delta^2$. Hence
$\chi_\pi^\sigma(\delta)=(2m+1)\sum_{\pi_H}\chi_{\pi_H}(\delta^2)$
at $\delta\in G_n^\sigma=H_n$.

We claim that for $\delta=\exp X$, $X\in\LieG_n^\sigma=\LieH_n$, we have 
$$\psi_E(\tr[X\upi^{-2n}\beta])=\psi(\tr[2X\upi^{-2n}\beta_H]).$$ 
For this we note that $\beta=\beta_H$ and $\psi_E(x)=\psi(x+\ov x)$.

Moreover we claim that 
$\upsi_n(\delta)=\upsi_{H,n}(\delta^2)$ for $\delta\in G_n^\sigma=H_n$.
For this we note that
$\upsi_n({}^tbu)=\psi_E((x+y)\upi^{-2n})$ if 
$u=\pmatrix 1&x&z\\0&1&y\\0&0&1\endpmatrix$, and this is 
$=\psi(2(x+\ov x)\upi^{-2n})$ if $y=\ov x$. But
$\psi_{H,n}({}^tbu)=\psi((x+\ov x)\upi^{-2n})$ at such $u\in U_H$ 
(thus with $y=\ov x$). 

Now $d(g^2)=dg$ when $p\not=2$. It follows that
$$1=\dim_{\C}\Hom_G(\ind_U^G\upsi,\pi)
=\int_{G_n^\sigma}\chi_\pi^\sigma(\delta)\psi_n(\delta)e_{G_n^\sigma}(\delta)$$
$$=\int_{H_n}(2m+1)\sum_{\pi_H}\chi_{\pi_H}(\delta^2)\psi_{H,n}(\delta^2)
e_{H_n}(\delta^2)$$
$$=(2m+1)\sum_{\pi_H}\dim_{\C}\Hom_H(\ind_{U_H}^H\upsi_H,\pi_H).$$
Hence $m=0$ and there is just one generic $\pi_H$ in the sum 
($\dim_{\C}\not=0$, necessarily $=1$).\enddemos

The excluded case of $p=2$ might follow on counting the factors of 2 in our 
argument.

We repeat the conclusion of [F]. Note that the unrestricted trace identity in 
proven in [F1], and the fundamental lemma in [F2]. Both proofs employ simple 
methods: the usage of regular-Iwahori functions in [F1] removes the need to 
compute and compare weighted orbital integrals, and explicit double coset
decomposition reduces the fundamental lemma in [F2] to elementary volume 
computations.

\proclaim{Corollary} Each discrete spectrum automorphic representation of
$\U(3,E/F)(\A)$ occurs in the discrete spectrum with multiplicity one. 
Each packet $($defined in $[\operatorname{F}])$ of such an infinite 
dimensional representation contains precisely one irreducible representation 
which is generic. Each packet $($defined in $[\operatorname{F}])$ of tempered 
admissible representations of $\U(3,E/F)$ contains precisely one irreducible 
representation which is generic. There are no generic representations in 
any $($nontempered$)$ quasipacket, locally and globally.\endproclaim

This type of argument, relating the number of Whittaker models to the 
values of the characters at suitable measures supported near the origin,
was first employed in [FK], which appeared in 1987, in an analogous 
situation of the metaplectic correspondence between $\GL(r)$ and its 
$n$-fold covering group.

\remark{Remark} On the space $L^2(\GL(n,F)\bs\GL(n,\A))$ of automorphic
forms $\phi$, we define the involution $\sigma$ by $(r(\sigma)\phi)(g)=
\phi(\sigma(g))$. On the space of Whittaker functions $W$ ($W(ng)=\psi(n)W(g)$,
$g\in\GL(n,\A)$, $n\in N(\A)$, where $N$ denotes the unipotent
upper triangular subgroup of $\GL(n)$), we choose the natural action of
$\sigma$, by ${}^\sigma W(g)=W(\sigma g)$. The map $\phi\mapsto W_\phi$,
where $W_\phi(g)=\int_{N(\A)\bs N(\A)}\phi(ng)\ov\psi(n)dn$, respects the
action of $\sigma$. 
Thus the global normalization of the action of $\sigma$ is the product of
the local normalizations ${}^\sigma W_v(g)=W_v(\sigma g)$.

Note also that underlying the character identity is an embedding of the 
local representation $\pi$, which is $\sigma$-elliptic, as a component of 
a global $\sigma$-invariant cuspidal representation $\Pi$ of $\GL(3,\A_E)$ 
where $F$ is a totally imaginary field. The character identity follows from 
using the trace formulae identity, applying ``generalized linear
independence of characters'' and thus isolating $\Pi$, and further, $\pi$.
The normalization of $\sigma$ on the generic ($\Pi$ and) $\pi$ implies
that no sign occurs in the character identity.

Had we not fixed the choice of the sign of $r(\sigma)$ at all places, the
twisted character in our trace identity (at our place) might in principle
be replaced by its negative. However, the identity
$$-\tr\pi(fdg\times\sigma)=(2m+1)\sum_{\pi_H}\tr\pi_H(f_Hdh),\qquad
-\chi_\pi^\sigma(\delta)=(2m+1)\sum_{\pi_H}\chi_{\pi_H}(\gamma),$$
cannot hold, as evaluating it with our 
$fdg=\psi_ne_{G_n^\sigma}$ and
$f_Hdh=\psi_{H,n}e_{H_n}$, $n$ large, would give $-1$ on the left,
and a nonnegative integer on the right. This provides an independent
verification that our normalization of the sign of $r(\sigma)$ is correct.

\endremark

\heading{3. Incomplete global proof}\endheading

The second proof of Proposition 3.5 of [F], on p. 48, is global, but
incomplete. The false assertion is on lines 21-22: ``Proposition 8.5(iii)
(p. 172) and 2.4(i) of [GP] imply that for some $\pi$ with $m(\pi)\not=0$
above, we have $m(\pi)=1$''. 
Indeed, [GP], Prop. 2.4, defines $L^2_{0,1}$ to be the orthocomplement
in the space $L^2_0$ (of cusp forms) of ``all hypercusp forms'', and
claims: ``(i) $L^2_{0,1}$ has multiplicity 1''. ([GP], 8.5 (iii), asserts
that $\pi$ is in $L^2_{0,1}$). Now the sentence of [F], p. 48, l. 21-22
assumes that [GP], 2.4(i), means that any irreducible $\pi$ in $L^2_{0,1}$ 
occurs in $L^2_0$ with multiplicity one. But the standard techniques of [GP], 
2.4, show only that any irreducible $\pi$ in $L^2_{0,1}$ occurs in $L^2_{0,1}$ 
with multiplicity one. {\it A-priori} there can exist $\pi'$ in $L^2_0$, 
isomorphic and orthogonal to $\pi\subset L^2_{0,1}$. In such a case we would 
have $m(\pi)>1$. Such a $\pi'$ is locally generic (all of its local components 
are generic), isomorphic to a generic cuspidal $\pi$, and the question 
boils down to whether this implies that $\pi'$ is generic (the linear form 
$L(\phi)=\int_{U_H(F)\bs U_H(\A)}\phi(u)\upsi_H(u)du$ is nonzero on 
$\pi\subset L_0^2$). This last claim might follow on using the theory 
of the Theta correspondence, but this has not been done as yet.
In summary, a clear form of [GP], 2.4(i) is: ``Any irreducible $\pi$ in 
$L^2_{0,1}$ occurs in $L^2_{0,1}$ with multiplicity one.'' 
In the analogous situation of GSp(2) such a statement is made in [So].
It is not sufficiently strong to be useful for us.

We noticed that the global argument of [F], p. 48, which was first
proposed in a preprint version of [F] in 1983, is incomplete while
generalizing it in [F3] to the context of the symplectic group, where
work of Kudla, Rallis, Langlands, Shahidi on the Siegel-Weil formula and
on L-functions is available to show that a locally generic cuspidal 
representation which is equivalent at almost all places to a generic 
cuspidal representation is generic. A local proof, based on a twisted 
analogue of Rodier's result, is also used in [F4], in the context of
the symmetric square lifting.

\heading{4. Review of Rodier's proof}\endheading

We shall reduce Theorem 2 to Theorem 1 for $G$ (not $H$), so we begin 
by recalling the main lines in Rodier's proof in the context of $G$.
Choose $d=\diag(\upi^{-r+1},\upi^{-r+3},\dots,\upi^{r-1})$. It lies
in the unitary group, namely $\sigma(d)=d$, since $\upi$ is in $F$. 
Put $V_n=d^nG_nd^{-n}$ and $\upsi_n(v)=\psi_n(d^{-n}vd^n)$ ($v\in V_n$). 
Recall that $\psi_n$ is defined to be supported on $G_n$.
Note that $\sigma(G_n)=G_n$, $\sigma(U_n)=U_n$, $\sigma\psi_n=\psi_n$, 
and that the entries in the $j$th line ($j\not=0$) above or below the 
diagonal of $v=(v_{ij})$ in $V_n$ lie in $\upi^{(1-2j)n}R_E$ (thus 
$v_{i,i+j}\in\upi^{(1-2j)n}R_E$ if $j>0$, and also when $j<0$). Thus 
$V_n\cap U$ is a $\sigma$-invariant strictly increasing sequence of 
compact and open subgroups of $U$ whose union is $U$, while 
$V_n\cap({}^tUH)$ -- where ${}^tUH$ is the lower triangular 
subgroup of $G$ -- is a strictly decreasing sequence of compact open 
subgroups of $G$ whose intersection is the element $I$ of $G$. 
Note that $\upsi_n=\upsi$ on $V_n\cap U$.

Consider the induced representations $\ind_{V_n}^G\upsi_n$, and the
intertwining operators
$$A_n^m:\ind_{V_n}^G\upsi_n\to\ind_{V_m}^G\upsi_m,$$
$$(A_n^m\f)(g)=((e_{V_m}\upsi_m)\ast \f)(g)
=\int_G\upsi_m(u)\f(u^{-1}g)e_{V_m}(u)$$
($g$ in $G$, $\f$ in $\ind_{V_n}^G\upsi_n$, $e_{V_m}=|V_m|^{-1}1_{V_m}dg$,
$|V_m|$ denotes the volume
of $V_m$, $1_{V_m}$ denotes the characteristic function of $V_m$).
For $m\ge n\ge 1$ we have
$$(A_n^m\f)(g)=((e_{V_m\cap U}\upsi)\ast \f)(g)
=\int_G\upsi(u)\f(u^{-1}g)e_{V_m\cap U}(u).$$
Hence $A_m^\ell\circ A_n^m=A_n^\ell$ for $\ell\ge m\ge n\ge 1$. So
$(\ind_{V_n}^G\upsi_n,A_n^m\,\,(m\ge n\ge 1))$ is an inductive system
of representations of $G$. Denote by $(I,A_n:\ind_{V_n}^G\upsi_n\to I)$
($n\ge 1$) its limit.

The intertwining operators $\phi_n:\ind_{V_n}^G\upsi_n\to\ind_U^G\upsi$,
$$(\phi_n(\f))(g)=(\upsi 1_U\ast \f)(g)=\int_U\upsi(u)\f(u^{-1}g)du,$$
satisfy $\phi_m\circ A_n^m=\phi_n$ if $m\ge n\ge 1$. Hence there exists
a unique intertwining operator $\phi:I\to\ind_U^G\upsi$ with $\phi\circ 
A_n=\phi_n$ for all $n\ge 1$. Proposition 3 of [R] asserts that

\proclaim{Lemma 1} The map $\phi$ is an isomorphism of $G$-modules.
\hfill$\square$\endproclaim

\proclaim{Lemma 2} There exists $n_0\ge 1$ such that 
$\upsi_n\ast\upsi_m\ast\upsi_n=|V_n||V_m\cap V_n|\upsi_n$ for all 
$m\ge n\ge n_0$.\endproclaim

\demo{Proof} This is Lemma 5 of [R]. We review its proof (the first displayed 
formula in the proof of this Lemma 5, [R], p. 159, line -8, should be erased).

There are finitely many representatives $u_i$ in $U\cap V_m$ for the
cosets of $V_m$ modulo $V_n\cap V_m$. Denote by $\varepsilon(g)$ the
Dirac measure in a point $g$ of $G$. Consider
$$(\varepsilon(u_i)\ast\upsi_n1_{V_m\cap V_n})(g)
=\int_G\varepsilon(u_i)(gh^{-1})(\upsi_n1_{V_m\cap V_n})(h)dh
=\upsi_n(u_i^{-1}g)=\upsi_m(u_i)^{-1}\upsi_m(g).$$
Note here that if the left side is nonzero, then $g\in u_i(V_m\cap V_n)
\subset V_m$. Conversely, if $g\in V_m$, then $g\in u_i(V_m\cap V_n)$
for some $i$. Hence $\upsi_m=\sum_i\upsi_m(u_i)\varepsilon(u_i)\ast
\upsi_n1_{V_m\cap V_n}$, thus
$$\upsi_n\ast\upsi_m\ast\upsi_n=\sum_i\upsi_m(u_i)\upsi_n\ast
\varepsilon(u_i)\ast\upsi_n1_{V_m\cap V_n}\ast\upsi_n.$$
Since $\upsi_n1_{V_m\cap V_n}\ast\upsi_n=|V_m\cap V_n|\upsi_n$, this is
$$=\sum_i\upsi_m(u_i)|V_m\cap V_n|\upsi_n\ast\varepsilon(u_i)\ast\upsi_n.$$
But the key Lemma 4 of [R] asserts that 
$\upsi_n\ast\varepsilon(u)\ast\upsi_n\not=0$ implies that $u\in V_n$.
Hence the last sum reduces to a single term, with $u_i=1$, and we obtain
$$=|V_m\cap V_n|\upsi_n\ast\upsi_n=|V_m\cap V_n||V_n|\upsi_n.$$

This completes the proof of the lemma.
\enddemos

\proclaim{Lemma 3} For an inductive system $\{I_n\}$ we have
$\Hom_G(\liminj I_n,\pi)=\limproj\Hom_G(I_n,\pi).$
\endproclaim

\demo{Proof} See, e.g., Rotman [Ro], Theorem 2.27. It is also verified in [R]. 
\enddemos

\proclaim{Corollary} We have 
$\quad\dim_\C\Hom_G(\ind_U^G\upsi,\pi)
=\lim_n|G_n|^{-1}\tr\pi(\psi_ndg).$
\endproclaim

\demo{Proof} As the $\dim_\C\Hom_G(\ind_{V_n}^G\upsi_n,\pi)$ are increasing 
with $n$, if they are bounded we get that they are independent of $n$ for 
sufficiently large $n$. Hence the left side of the corollary is equal to
$\lim_n\dim_\C\Hom_G(\ind_{V_n}^G\upsi_n,\pi)$. This is equal to
$\lim_n\dim_\C\Hom_G(\ind_{G_n}^G\psi_n,\pi)$ since 
$\upsi_n(v)=\psi_n(d^{-n}vd^n)$. This is equal to
$\lim_n\dim_\C\Hom_{G_n}(\psi_n,\pi|G_n)$ by Frobenius reciprocity.
This is equal to the right side of the corollary since 
$|G_n|^{-1}\pi(\psi_n dg)$ is a projection from $\pi$ to the space
of $\xi$ in $\pi$ with $\pi(g)\xi=\psi_n(g)\xi$ ($g\in G_n$), a space whose 
dimension is then $|G_n|^{-1}\tr\pi(\psi_n dg)$.
\enddemos

\heading{5. The twisted case}\endheading

We now reduce Theorem 2 to Theorem 1 for $G$.
Note that since $\sigma\upsi_n=\upsi_n$, the representations 
$\ind_{V_n}^G\upsi_n$ are $\sigma$-invariant, where $\sigma$ acts
on $\f\in\ind_{V_n}^G\upsi_n$ by $\f\mapsto\sigma \f$, $(\sigma \f)(g)=
\f(\sigma g)$. Similarly $\sigma\upsi=\upsi$ and $\ind_U^G\upsi$ is
$\sigma$-invariant. We then extend these representations ind of $G$
to the semidirect product $G'=G\rtimes\langle\sigma\rangle$ by putting
$(i(\sigma)\f)(g)=\f(\sigma(g))$. 

Let $\pi$ be an irreducible admissible representation of $G$ which is
$\sigma$-invariant. Thus there exists an intertwining operator
$A:\pi\to{}^\sigma\pi$, where ${}^\sigma\pi(g)=\pi(\sigma(g))$, with
$A\pi(g)=\pi(\sigma(g))A$. Then $A^2$ commutes with every $\pi(g)$ 
($g\in G$), hence $A^2$ is a scalar by Schur's lemma, and can be
normalized to be 1. This determines $A$ up to a sign. We extend $\pi$ 
from $G$ to $G'=G\rtimes\langle\sigma\rangle$ by putting $\pi(\sigma)=A$
once $A$ is chosen.

If $\Hom_G(\ind_U^G\upsi,\pi)\not=0$, its dimension is 1. Choose a 
generator $\ell:\ind_U^G\upsi\to\pi$. Define $A:\pi\to\pi$ by $A\ell(\f)
=\ell(i(\sigma)\f)$. Then
$$\Hom_G(\ind_U^G\upsi,\pi)=\Hom_{G'}(\ind_U^G\upsi,\pi).$$

Similarly we have
$\quad\Hom_G(\ind_{V_n}^G\upsi_n,\pi)=\Hom_{G'}(\ind_{V_n}^G\upsi_n,\pi).$

The right side in the last equality can be expressed as
$$\Hom_{G'}(\ind_{G_n}^G\psi_n,\pi)=\Hom_{G'_n}(\psi'_n,\pi|G'_n)\qquad
(G'_n=G_n\rtimes\langle\sigma\rangle).$$
The last equality follows from Frobenius reciprocity, where we extended
$\psi_n$ to a character $\psi'_n$ on $G'_n$ by $\psi'_n(\sigma)=1$. Thus 
$\psi'_n=\psi^1_n+\psi_n^\sigma$, with $\psi_n^\alpha(g\times\beta)=
\delta_{\alpha\beta}\psi_n(g)$, $\alpha$, $\beta\in\{1,\sigma\}$.

In this case $\Hom_{G'_n}(\psi'_n,\pi|G'_n)$ is isomorphic to the space 
$\pi_1$ of vectors $\xi$ in $\pi$ with $\pi(g)\xi=\psi_n(g)\xi$ for all 
$g$ in $G'_n$. In particular $\pi(g)\xi=\psi_n(g)\xi$ for all $g$ in $G_n$, 
and $\pi(\sigma)\xi=\xi$. Clearly $|G'_n|^{-1}\pi(\psi'_ndg')$ is a 
projection from the space of $\pi$ to $\pi_1$ (it is independent of the 
choice of the measure $dg'$). Its trace is then the dimension of the space 
Hom. We conclude a twisted analogue of the theorem of [R]:

\proclaim{Proposition 1} The integer $\dim_\C\Hom_{G'}(\ind_U^G\upsi,\pi)$
is equal to $|G'_n|^{-1}\tr\pi(\psi'_ndg')$ for all sufficiently large $n$.
\hfill$\square$\endproclaim

Note that $G'_n$ is the semidirect product of $G_n$ and the two-element
group $\langle\sigma\rangle$. With the natural measure assigning 1 to each
element of the discrete group $\langle\sigma\rangle$, we have $|G'_n|
=2|G_n|$. The result is then, for all sufficiently large $n$,
$${1\over 2}\tr\pi(\psi_ne_{G_n})
+{1\over 2}\tr\pi(\psi_ne_{G_n}\times\sigma)$$
(as $\psi'_n=\psi^1_n+\psi_n^\sigma$, $\psi^1_n=\psi_n$ and 
$\tr\pi(\psi_n^\sigma dg)=\tr\pi(\psi_ndg\times\sigma)$). By (the nontwisted)
Rodier's Theorem 1,
$$\dim_\C\Hom_G(\ind_U^G\upsi,\pi)=\lim_n\tr\pi(\psi_ne_{G_n}),$$
we conclude 

\proclaim{Proposition 2} We have
$\quad\dim_\C\Hom_{G'}(\ind_U^G\upsi,\pi)
=\lim_n\tr\pi(\psi_ne_{G_n}\times\sigma)\quad$ 
for all $\sigma$-invariant irreducible representations $\pi$ of $G$.
\hfill$\square$\endproclaim

The terms in the limit on the right can be written in terms of 
Harish-Chandra's twisted character, as
$$\int_G\chi_\pi^\sigma(g)\psi_n(g)e_{G_n}(g).$$
Again, put $e_{G_n^\sigma}=|G_n^\sigma|^{-1}\ch_{G_n^\sigma}dg$ where 
$\ch_{G_n^\sigma}$ is the characteristic function of $G_n^\sigma$ in $G$.

\proclaim{Proposition 3} The last displayed integral is equal to
$$\int_{G_n^\sigma}\chi_\pi^\sigma(g)\psi_n(g)e_{G^\sigma_n}(g).$$
\endproclaim

\demo{Proof} Consider the map $G_n^\sigma\times G_n^\sigma\bs G_n\to G_n$, 
$(u,k)\mapsto k^{-1}u\sigma(k)$. It is a closed immersion. More generally, 
given a semisimple element $s$ in a group $G$, we can consider the map 
$Z_{G^0}(s)\times Z_{G^0}(s)\bs G^0\to G^0$ by $(u,k)\mapsto k^{-1}usks^{-1}$. 
Our example is: $(s,G)=(\sigma,G_n\times\langle\sigma\rangle)$. 

Our map is in fact an analytic isomorphism since $G_n$ is a small 
neighborhood of the origin, where the exponential $e:\LieG_n\to G_n$ is an 
isomorphism. Indeed, we can transport the situation to the Lie algebra 
$\LieG_n$. Thus we write $k=e^Y$, $u=e^X$, $\sigma(k)=e^{(d\sigma)(Y)}$,
$k^{-1}u\sigma(k)=e^{X-Y+(d\sigma)(Y)}$, up to smaller terms. Here
$(d\sigma)(Y)=-J^{-1}{}^t\ov{Y}J$. So we just need to show that 
$(X,Y) \mapsto X-Y+(d\sigma)(Y)$, $Z_{\LieG_n}(\sigma)+
\LieG_n(\mod Z_{\LieG_n}(\sigma))\to\LieG_n$, is bijective. But this is 
obvious since the kernel of $(1-d\sigma)$ on $\LieG_n$ is precisely 
$Z_{\LieG_n}(\sigma)=\{X\in\LieG_n;(d\sigma)(X)=X\}$.

Changing variables on the terms 
on the right of Proposition 2 we get the equality:
$$\int_{G_n}\chi_\pi^\sigma(g)\psi_n(g)e_{G_n}(g)$$
$$=|G_n|^{-1}\int_{G_n^\sigma}\int_{G_n^\sigma\bs G_n}
\chi_\pi^\sigma(k^{-1}u\sigma(k))\psi_n(k^{-1}u\sigma(k))dkdu.$$
But $\sigma\psi_n=\psi_n$, $\psi_n$ is a homomorphism (on $G_n$), $G_n$ is 
compact, and $\chi_\pi^\sigma$ is a $\sigma$-conjugacy class function, so we 
end up with the expression of the proposition. Note that $\chi_\pi^\sigma$
is locally integrable on $G_n^\sigma$ and locally constant on its regular
set by the character relation stated in the proof of Prop. 3.5 of [F] above.
The proposition, and Theorem 2, follow.
\enddemos

\heading{6. Appendix. Germs of twisted characters}\endheading

Harish-Chandra [HC] showed that $\chi_\pi$ is locally integrable (Thm 1,
p. 1) and has a germ expansion near each semisimple element $\gamma$
(Thm 5, p. 3), of the form: 
$$\chi_\pi(\gamma\exp X)
=\sum_{\Cal O} c_\gamma({\Cal O},\pi)\wh\mu_{\Cal O}(X).$$
Here ${\Cal O}$ ranges over the nilpotent orbits in the Lie algebra $\LieM$
of the centralizer $M$ of $\gamma$ in $G$, $\mu_{\Cal O}$ is an invariant 
distribution supported on the orbit ${\Cal O}$, 
$\wh\mu_{\Cal O}$ is its Fourier 
transform with respect to a symmetric nondegenerate $G$-invariant bilinear 
form $B$ on $\LieM$ and a selfdual measure, and $c_\gamma({\Cal O},\pi)$ are 
complex numbers. Both $\mu_{\Cal O}$ and 
$c_\gamma({\Cal O},\pi)$ depend on a choice of
a Haar measure $d_{\Cal O}$ on the centralizer 
$Z_G(X_0)$ of $X_0\in{\Cal O}$, but their product does not. 
The $X$ ranges over a small neighborhood of the origin in $\LieM$. We shall
be interested only in the case of $\gamma=1$, and thus omit $\gamma$ from
the notations. The size of the domain where the germ expansion holds is 
studied in Waldspurger [W]. 

Suppose that $G$ is quasisplit over $F$, and $U$ is the unipotent radical
of a Borel subgroup $B$. Let $\upsi:U\to\C^1$ be the nondegenerate
character of $U$ (its restriction to each simple root subgroup is 
nontrivial) specified in Rodier [R], p. 153. The number 
$\dim_{\C}\Hom(\ind_U^G\upsi,\pi)$ of $\upsi$-Whittaker functionals on 
$\pi$ is known to be zero or one. Let $\LieG_0$ be a selfdual lattice 
in the Lie algebra $\LieG$ of $G$. Denote by $\ch_0$ the characteristic
function of $\LieG_0$ in $\LieG$. Rodier [R], p. 163, showed that there 
is a regular nilpotent orbit ${\Cal O}={\Cal O}_\upsi$ 
such that $c({\Cal O},\pi)$ is
not zero iff $\dim_{\C}\Hom(\ind_U^G\upsi,\pi)$ is one, in fact
$\wh\mu_{\Cal O}(\ch_0)c({\Cal O},\pi)$ is one in this case. 
Alternatively put, normalizing $\mu_{\Cal O}$ by 
$\wh\mu_{\Cal O}(\ch_0)=1$, we have $c({\Cal O},\pi)=
\dim_{\C}\Hom(\ind_U^G\upsi,\pi)$. This is shown in [R] for all $p$ if 
$G=\GL(r,F)$, and for general quasisplit $G$ for all $p\ge 1+
2\sum_{\alpha\in S}n_\alpha$, if the longest root is 
$\sum_{\alpha\in S}n_\alpha\alpha$ in a basis $S$ of the root system. 
A generalization of Rodier's theorem to degenerate Whittaker models 
and nonregular nilpotent orbits is given in Moeglin-Waldspurger [MW]. 
See [MW], I.8, for the normalization of measures. In particular they 
show that $c({\Cal O},\pi)> 0$ for the nilpotent orbits 
${\Cal O}$ of maximal dimension with $c({\Cal O},\pi)\not=0$. 
For applications to minimal representations see Savin [S].

Harish-Chandra's results extend to the twisted case. The twisted character 
is locally integrable (Clozel [C], Thm 1, p. 153), and there exist unique 
complex numbers $c^\theta({\Cal O},\pi)$ ([C], Thm 3, p. 154) with
$\chi_\pi^\theta(\exp X)=\sum_{\Cal O} 
c^\theta({\Cal O},\pi)\wh\mu_{\Cal O}(X).$
Here ${\Cal O}$ ranges over the nilpotent orbits in the 
Lie algebra $\LieG^\theta$ of the group $G^\theta$ of the 
$g\in G$ with $g=\theta(g)$. Further, $\mu_{\Cal O}$ is an invariant 
distribution supported on the orbit ${\Cal O}$ (it is unique 
up to a constant, not unique as stated in [HC], Thm 5, and [C], Thm 3); 
$\wh\mu_{\Cal O}$ is its Fourier transform, and $X$ ranges over a small 
neighborhood of the origin in $\LieG^\theta$. 

In this section we compute the expression displayed in Proposition 
3 using the germ expansion $\chi_\pi^\sigma(\exp X)=\sum_{\Cal O} 
c^\sigma({\Cal O},\pi)\wh\mu_{\Cal O}(X).$ This expansion 
means that for any test measure $fdg$ supported on a 
small enough neighborhood of the identity in $G$ we have
$$\int_{\LieG^\sigma}f(\exp X)\chi_\pi^\sigma(\exp X)dX
=\sum_{\Cal O}c^\sigma({\Cal O},\pi)\int_{\Cal O}
[\int_{\LieG^\sigma}f(\exp X)\psi(\tr(XZ))dX] d\mu_{\Cal O}(Z).$$
Here ${\Cal O}$ ranges over the nilpotent orbits in $\LieG^\sigma$, 
$\mu_{\Cal O}$ is an invariant distribution supported on the orbit 
${\Cal O}$, $\wh\mu_{\Cal O}$ is its Fourier transform. 
The $X$ range over a small neighborhood of the origin in $\LieG^\sigma$.
Since we are interested only in the case of the unitary group, and to 
simplify the exposition, we take $G=\GL(r,E)$ and the involution $\sigma$ 
whose group of fixed points is the unitary group $H=\U(r,E/F)$. 
In this case there is a unique regular nilpotent orbit ${\Cal O}_0$.

We normalize the measure $\mu_{\Cal O_0}$ on the orbit ${\Cal O}_0$ of $\beta$ 
in $\LieG^\sigma$ by the requirement that $\wh\mu_{\Cal O_0}(\ch_0^\sigma)$ 
is 1, thus that $\int_{\beta+\upi^n\LieG_0^\sigma}d\mu_{\Cal O_0}(X)
=q^{n\dim(\Cal O_0)}$ for large $n$. Equivalently a measure on
an orbit $\Cal O\simeq G/Z_G(Y)$ ($Y\in\Cal O$) is defined by a
measure on its tangent space $m=\LieG/Z_{\LieG}(Y)$ ([MW], p. 430)
at $Y$, taken to be the selfdual measure with respect to the symmetric 
bilinear nondegenerate $F$-valued form $B_Y(X,Z)=\tr(Y[X,Z])$ on $m$.

\proclaim{Proposition 4} If $\pi$ is a $\sigma$-invariant admissible
irreducible representation of $G$ and ${\Cal O}_0$ is the regular nilpotent
orbit in $\LieG^\sigma$, then the coefficient $c^\sigma({\Cal O}_0,\pi)$ 
in the germ expansion of the $\sigma$-twisted character $\chi_\pi^\sigma$ 
of $\pi$ is equal to
$$\dim_\C\Hom_{G'}(\ind_U^G\upsi,\pi)
=\dim_\C\Hom_{G}(\ind_U^G\upsi,\pi).$$
This number is one if $\pi$ is generic, and zero otherwise.
\endproclaim

\demo{Proof} We compute the expression displayed in Proposition 3 
as in [MW], I.12. It is a sum over the nilpotent orbits $\Cal O$ 
in $\LieG^\sigma$, of $c^\sigma(\Cal O,\pi)$ times
$$|G_n^\sigma|^{-1}\wh\mu_{\Cal O}(\psi_n\circ e)
=|G_n^\sigma|^{-1}\mu_{\Cal O}(\wh{\psi_n\circ e})
=|G_n^\sigma|^{-1}\int_{\Cal O}\wh{\psi_n\circ e}(X)d\mu_{\Cal O}(X).$$
The Fourier transform (with respect to the character $\psi_E$) of 
$\psi_n\circ e$,
$$\wh{\psi_n\circ e}(Y)=\int_{\LieG^\sigma}\psi_n(\exp Z)\ov\psi_E(\tr ZY)dZ
=\int_{\LieG_n^\sigma}\psi_E(\tr Z(\upi^{-2n}\beta-Y))dZ,$$
is the characteristic function of $\upi^{-2n}\beta+\upi^{-n}\LieG_0^\sigma
=\upi^{-2n}(\beta+\upi^n\LieG_0^\sigma)$ multiplied by the volume 
$|\LieG_n^\sigma|=|G_n^\sigma|$ of $\LieG_n^\sigma$. Hence we get
$$=\int_{\Cal O\cap(\upi^{-2n}(\beta+\upi^n\LieG_0^\sigma))}d\mu_{\Cal O}(X)
=q^{n\dim(\Cal O)}\int_{\Cal O\cap(\beta+\upi^n\LieG_0^\sigma)}
d\mu_{\Cal O}(X).$$
The last equality follows from the homogeneity result of [HC], Lemma 3.2, 
p. 18. For sufficiently large $n$ we have that
$\beta+\upi^n\LieG_0^\sigma$ is contained only in the orbit $\Cal O_0$ of
$\beta$. Then only the term indexed by $\Cal O_0$ remains in the
sum over $\Cal O$, and 
$$\int_{\Cal O_0\cap(\beta+\upi^n\LieG_0^\sigma)}d\mu_{\Cal O_0}(X)
=\int_{\beta+\upi^n\LieG_0^\sigma}d\mu_{\Cal O_0}(X)$$
equals $q^{-n\dim(\Cal O_0)}$ (cf. [MW], end of proof of Lemme I.12).
The proposition follows.\enddemos
\bigskip
\def\refe#1#2{\n\hangindent 5em\hangafter1\hbox to 5em{\hfil#1\quad}#2}
\subheading{References}
\medskip

\refe{[BZ1]}{J. Bernstein, A. Zelevinskii, Representations of the
group $\GL(n,F)$ where $F$ is a nonarchimedean local field, {\it 
Russian Math. Surveys} 31 (1976), 1-68.}

\refe{[BZ2]}{J. Bernstein, A. Zelevinsky, Induced representations of 
reductive $p$-adic groups, {\it Ann. Sci. \'Ecole Norm. Sup.} 10 (1977), 
441-472.}

\refe{[C]}{L. Clozel, Characters of non-connected, reductive $p$-adic
groups, {\it Canad. J. Math.} 39 (1987), 149-167.}


\refe{[F]}{Y. Flicker, Packets and liftings for U(3), {\it J. Analyse Math.} 
50 (1988), 19-63.}

\refe{[F1]}{Y. Flicker, Base change trace identity for U(3), {\it J. Analyse 
Math.} 52 (1989), 39-52.}

\refe{[F2]}{Y. Flicker, Elementary proof of a fundamental lemma for a 
unitary group, {\it Canad. J. Math.} 50 (1998), 74-98.}

\refe{[F3]}{Y. Flicker, {\it Lifting Automorphic Forms on} PGSp(2) {\it and} 
SO(4) {\it to} PGL(4), research monograph, 2001; see also: Automorphic forms 
on PGSp(2), {\it Elect. Res. Ann. Am. Math. Soc.} 10 (2004), 39-50.
http://www.ams.org/era/}

\refe{[F4]}{Y. Flicker, {\it Automorphic Representations of Low Rank Groups},
research monograph, 2003.}

\refe{[FK]}{Y. Flicker, D. Kazhdan, Metaplectic correspondence,
{\it Publ. Math. IHES} 64 (1987), 53-110}

\refe{[GP]}{S. Gelbart, I. Piatetski-Shapiro, Automorphic forms and 
L-functions for the unitary group, in {\it Lie Groups Representations II},
Springer Lecture Notes 1041 (1984), 141-184.}



\refe{[H]}{G. Harder, {\it Eisensteinkohomologie und die Konstruktion
gemischter Motive},\hb Springer Lecture Notes 1562 (1993).}

\refe{[HC]}{Harish-Chandra, {\it Admissible invariant distributions on
reductive $p$-adic groups}, notes by S. DeBacker and P. Sally, AMS Univ.
Lecture Series 16 (1999); see also: {\it Queen's Papers in Pure and Appl. 
Math.} 48 (1978), 281-346.}



\refe{[MW]}{C. Moeglin, J.L. Waldspurger, Mod\`eles de Whittaker 
d\'eg\'en\'er\'es pour des groupes $p$-adiques, {\it Math. Z.} 196
(1987), 427-452.}


\refe{[R]}{F. Rodier, Mod\`ele de Whittaker et caract\`eres de 
repr\'esentations, {\it Non-commuta-tive harmonic analysis}, Springer
Lecture Notes 466 (1975), 151-171.}


\refe{[Ro]}{J. Rotman, {\it An Introduction to Homological Algebra},
Academic Press 1979.}

\refe{[S]}{G. Savin,  Dual pair $G_J\times\PGL_2$ where $G_J$ is the 
automorphism group of the Jordan algebra $J$, {\it  Invent. Math.} 118 (1994),
141-160.}

\refe{[So]}{D. Soudry, A uniqueness theorem for representations of GSO(6) 
and the strong multiplicity one theorem for generic representations of GSp(4),
{\it Israel J. Math.} 58 (1987), 257-287.} 

\refe{[W]}{J.-L. Waldspurger, Homog\'en\'eit\'e de certaines distributions 
sur les groupes $p$-adiques, {\it Publ. Math. IHES} 81 (1995), 25-72.} 

\enddocument